# An inverse problem formulation for parameter estimation of a reaction-diffusion model of low grade gliomas

**Amir Gholami · Andreas Mang · George Biros**

**Abstract** We present a numerical scheme for solving a parameter estimation problem for a model of low-grade glioma growth. Our goal is to estimate the spatial distribution of tumor concentration, as well as the magnitude of anisotropic tumor diffusion. We use a constrained optimization formulation with a reaction-diffusion model that results in a system of nonlinear partial differential equations (PDEs). In our formulation, we estimate the parameters using partially observed, noisy tumor concentration data at two different time instances, along with white matter fiber directions derived from diffusion tensor imaging (DTI). The optimization problem is solved with a Gauss-Newton reduced space algorithm. We present the formulation and outline the numerical algorithms for solving the resulting equations. We test the method using a synthetic dataset and compute the reconstruction error for different noise levels and detection thresholds for monofocal and multifocal test cases.

## 1 Introduction

Gliomas are tumors that arise from glial cells in the brain. They account for 29% of all brain and central nervous system (CNS) tumors, and 80% of all malignant brain tumors of about 60,000 cases diagnosed each year in the United States (Dolecek et al., 2012). Despite advances in surgery, chemotherapy, and radiotherapy, the median survival rate with therapy has remained about one

Institute for Computational Engineering and Sciences, The University of Texas at Austin, TX 78712, USA
E-mail: i.amirgh@gmail.com; andreas@ices.utexas.edu; biros@ices.utexas.edu

year in the past 30 years (Newton, 1994; Salcman, 1980; Seither et al., 1995; Swanson et al., 2008; Wrensch et al., 2002).

One of the key challenges in treating gliomas is their aggressive infiltration into the healthy tissue, well beyond the visible bulk of the tumor in standard clinical imaging modalities (Holland et al., 1985; Seither et al., 1995). Thus it is hard to decide on how much tissue to resect in surgery or radiate in radiotherapy. Considering a large margin in radiotherapy may destroy healthy tissue, while small margins may result in faster recurrence of the tumor (Nazzaro and Neuwelt, 1990; Silbergeld and Chicoine, 1997). In clinical practice, a margin of 1cm–4cm is typically used (Hochberg and Pruitt, 1980; Lawrence et al., 2010). This margin is derived by statistical analysis of clinical data (Hochberg and Pruitt, 1980). But there is no consensus about its particular value and there is no systematic way to select this margin for an individual patient.

Our work has the potential to provide guidance on fine tuning the treatment margin in a systematic way, which in turn may improve treatment. It is based on integrating mathematical models of tumor growth with imaging, using automatic parameter calibration.

*Contributions.* We address the problem of determining the full extent of the tumor infiltration for a reaction-diffusion tumor growth model. We summarize our contributions below.

- We present the mathematical formulation of PDE-constrained parameter estimation for tumor growth. We invert for the initial condition of the tumor concentration and the extent of anisotropic infiltration. Unlike existing approaches, our parameterization allows for calibration of multifocal gliomas.
- We present and verify a numerical scheme for the solution of the parameter estimation problem. Due to the 3D PDEs that describe the tumor growth the estimation problem is challenging. We present novel numerical schemes that enable fast solution of this problem.
- We apply our scheme to estimate of the extent of tumor infiltration. We account for the fact that we have only limited information on the actual extent of the tumor in the patient data. We report preliminary results for synthetic tumors, and the associated reconstruction errors.

Including estimation of the anisotropic infiltration rate makes the problem more complex, but it is inevitable as there is currently no patient specific method to measure it[1]. In our experiments the data consists of two or more noisy (in terms of errors), partially observed, segmented images of the tumor. The partial observation corresponds to the visible bulk of the tumor observed in the images.

*Limitations.* The primary focus of this paper is on the formulation, numerical analysis, and feasibility of the proposed method. We use a widely adopted reaction-diffusion model for tumor growth (Jbabdi et al., 2005; Swanson et al., 2000, 2008, 2002). We do not consider mass effect (deformations of the parenchyma due to tumor growth). The use of more complex tumor growth models accounting for mass effect, edema, necrosis, angiogenesis and chemotaxis remains to be investigated (Engwer et al., 2014; Habib et al., 2003; Hawkins-Daarud et al., 2013; Hogea et al., 2008b). Moreover our method requires at least two consecutive time frames of tumor growth. This is a limitation especially for high grade tumors since there is rarely two time frames available. However, treatment is sometimes delayed for low grade tumors. These tumors are left in an observed state, for which multiple time frames may be available. Also, we require actual tumor concentration values. In practice one can get tumor classes from segmentation but not the concentration values (Mang et al., 2012). However, some methods of approximating the concentration values from ADC (apparent diffusion coefficient) data have been proposed (Anderson et al., 2000; Atuegwu et al., 2013; Weis et al., 2013).

[1] Manual tuning for this value has been reported (Jbabdi et al., 2005; Sodt et al., 2014).

*Related Work.* There is a long tradition in the design of mathematical models for cancer progression. The complexity of the underlying bio-physiology results in a diversity of mathematical models, accounting for phenomena on the molecular, cellular and/or tissue scale (Bellomo et al., 2008). We limit this review to tissue scale methods for gliomas and their integration with clinical images.

Most work attempting to link mathematical models to images is based on reaction-diffusion type equations (Murray, 1989)(see also §2). Its main assumption is that cancerous cells grow and infiltrate tissue due to cell division (proliferation) and migration (which can be modeled by diffusion). Although this model is simplistic and purely phenomenological, it has been shown that it can capture tumor progression on a tissue level (Clatz et al., 2005; Konukoglu et al., 2010a,b; Mang et al., 2012; Rekik et al., 2013).

This model has been used extensively for recovering tumor growth patterns in individual patients (Clatz et al., 2005; Hogea et al., 2007, 2008b; Jbabdi et al., 2005; Konukoglu et al., 2010b; Mang et al., 2012; Rekik et al., 2013; Swanson et al., 2000, 2002), for estimating the physiological (true) tumor boundary (Cobzas et al., 2009; Konukoglu et al., 2010a; Mosayebi et al., 2012), and for studying effects of clinical intervention (Powathil et al., 2007; Rockne et al., 2009, 2010; Swanson et al., 2008; Tracqui et al., 1995). The parameter calibration is typically driven using segmented data of patient images. Besides only considering tissue composition (Hogea et al., 2007, 2008b; Swanson et al., 2000, 2002), white matter architecture (neuronal pathways), obtained from diffusion tensor imaging (DTI) data, has also been included into the models, in an attempt to account for the experimentally observed prevalent migration of cancerous cells along the white matter tracts (Bondiau et al., 2008; Clatz et al., 2005; Cobzas et al., 2009; Engwer et al., 2014; Konukoglu et al., 2010a,b; Mang et al., 2012; Mosayebi et al., 2012; Painter and Hillen, 2013; Rekik et al., 2013). Models that account for the mechanical interaction of the tumor with its surroundings have been described in (Clatz et al., 2005; Hogea et al., 2007; Mohamed and Davatzikos, 2005).

A common approach for parameter estimation is manual calibration (Clatz et al., 2005; Swanson et al., 2002). However, it is somewhat difficult to reproduce results and this approach does not scale with the number of parameters. For this reason, several groups have developed automatic parameter calibration algorithms. One approach is to use an asymptotic approximation of the reaction diffusion type equations on the basis of a traveling wave solution (Swanson et al., 2008). The associated analytical result establishes a connection between the velocity of the (spherical) tumor front and the parameters of the underlying reaction diffusion type model. An extension that accounts for the heterogeneity of the tissue as well as the structure of white matter pathways has been described in (Konukoglu et al., 2010a,b). This method has been applied to time series of images (Konukoglu et al., 2010b) and to imaging studies based on a single snapshot in time (Rekik et al., 2013).

We use the general framework of parameter estimation, based on optimal control theory (Biros and Ghattas, 2005a; Biros and Ghattas, 2005b; Hogea et al., 2008b; Bangerth and Joshi, 2008). Our approach extends the work (Hogea et al., 2008b) in that we use second order (Hessian) information for the numerical optimization, instead of using first order information only or falling back to a derivative free optimization (Mang et al., 2012). We use a series of novel preconditioners which noticeably speeds up the time to solution. Also the experiments of (Hogea et al., 2008b) were limited to the 1D case, only. We report results for the 2D and 3D case, in addition to accounting for anisotropic diffusion in white matter.

*Outline of the paper.* In section §2, we introduce the tumor growth model (i.e. the forward problem). In §3, we present the optimization problem to solve for the unknowns (i.e. the inverse problem), and in §4 we describe details of the related numerical methods. Finally, in §5 we test the proposed algorithm accounting for different noise levels and different detection thresholds. We conclude with a discussion of the designed method and the reported results in §6.

## 2 Tumor Model

The tumor growth model we use has been widely adopted in the literature to model the spatio-temporal spread of cancerous cells on a tissue level (Hogea et al., 2007; Jbabdi et al., 2005; Murray, 1989; Swanson et al., 2000, 2008; Tracqui et al., 1995). It can be stated as follows:

Rate of change of cells in time = proliferation rate + motility (diffusion) rate.

The equivalent mathematical formulation is given by the following partial differential equation:

$$\frac{\partial c}{\partial t} - Dc - R(c) = 0 \text{ in U}, \tag{1}$$

$$\frac{\partial c}{\partial n} = 0 \text{ on } \Gamma \times (0{,}1), \tag{2}$$

subjected to an initial distribution $c(t = 0) = c_0$. Here, $c$ is the normalized tumor concentration (i.e. $c \in (0,\ 1])$, $D$ is a linear differential operator modeling the migration of the tumor cells, and $R(c)$ is a nonlinear reaction term, modeling proliferation and necrosis of the tumor cells. Furthermore, $U := \mathcal{B} \times (0, 1]$, where $\mathcal{B}$ is the spatial domain of brain and $(0, 1]$ is the non-dimensional time interval. Moreover, $\Gamma$ refers to the boundaries of CSF and skull of the brain where the tumor cells do not infiltrate into (Tracqui et al., 1995).

The differential operator $D$ that models tumor infiltration is based on a model of inhomogeneous, anisotropic diffusion:

$$Dc = \nabla \cdot (\mathbf{K}(\mathbf{x})\nabla c), \tag{3}$$

where

$$\mathbf{K}(\mathbf{x}) = k_0(\mathbf{x})\mathbf{I} + k_f\mathbf{T}(\mathbf{x}). \tag{4}$$

Here, $k_0(\mathbf{x})$ captures the inhomogeneity due to different diffusion rates in white and grey matter, and $\mathbf{x} = (x, y, z) \in \mathcal{B}$. The inhomogeneity of tumor infiltration follows from the experimental study in (Giese et al., 1996), in which it was observed that glioma cells have a higher motility rate in white matter as compared to grey matter. This observation has recently been related to the higher cell density in grey matter (Sodt et al., 2014). To account for these different motility rates, the diffusion coefficient is assumed to be inhomogeneous (i.e. $k_0$ is a function of the location $\mathbf{x}$), with a higher diffusion rate in white matter than in grey matter (Silbergeld and Chicoine, 1997; Swanson et al., 2000, 2002; Tracqui et al., 1995). Therefore, we have

$$k_0(\mathbf{x}) = \begin{cases} k_w, & \text{if } \mathbf{x} \text{ in white matter} \\ k_g, & \text{if } \mathbf{x} \text{ in grey matter} \\ 0, & \text{otherwise} \end{cases}$$

We use a five fold difference, i.e. $\frac{k_w}{k_g} = 5$ (Hogea et al., 2008a; Swanson et al., 2002).

Similarly, it has been suggested that glioma cells have a directional preference in their infiltration. This can be accounted for by introducing an anisotropic diffusion operator (Horsfield and Jones, 2002; Jbabdi et al., 2005; Painter and Hillen, 2013; Stadlbauer et al., 2009). The second term in Eq. 4, $\mathbf{T}(\mathbf{x})$, captures this behavior. $\mathbf{T}(\mathbf{x})$ is the weighted diffusion tensor derived from diffusion tensor imaging (DTI) data[2]. We compute $\mathbf{T}(\mathbf{x})$, by scaling the eigendirections and eigenvalues derived from the $3 \times 3$ DTI tensor by the so called fractional anisotropy ($FA$):

$$\mathbf{T} = FA(\lambda_1\mathbf{e_1}\mathbf{e_1}^T + \lambda_2\mathbf{e_2}\mathbf{e_2}^T + \lambda_3\mathbf{e_3}\mathbf{e_3}^T), \tag{5}$$

[2] DTI is an MR imaging technique that measures water diffusion tensor at every point in the brain (Le Bihan et al., 2001)

where $FA$ is given by

$$FA = \sqrt{\frac{1}{2}} \frac{\sqrt{(\lambda_1 - \lambda_2)^2 + (\lambda_2 - \lambda_3)^2 + (\lambda_3 - \lambda_1)^2}}{\sqrt{\lambda_1^2 + \lambda_2^2 + \lambda_3^2}}. \tag{6}$$

Here, $\mathbf{e_1}, \mathbf{e_2}, \mathbf{e_3}$ and $\lambda_1, \lambda_2, \lambda_3$ are the corresponding eigendirections and eigenvalues, computed at each point in the brain. The diffusion tensor $\mathbf{T}$ is additionally scaled by $k_f$ to account for differences in the diffusion rates between tumor cells and water. Other methods of deriving $\mathbf{T}$ by scaling only the principle direction have been suggested as well (Jbabdi et al., 2005; Painter and Hillen, 2013; Sodt et al., 2014). For example, $\mathbf{T}$ can have the form of:

$$\mathbf{T} = \lambda_1 \mathbf{e_1} \mathbf{e_1}^T. \tag{7}$$

In this approach only the principle direction is preferred when there is fiber crossings (the planar case). To adress this, (Jbabdi et al., 2005) suggested a linear scaling of the form:

$$\mathbf{T} = (a_1 \mathbf{e_1} \mathbf{e_1}^T + a_2 \mathbf{e_2} \mathbf{e_2}^T + a_3 \mathbf{e_3} \mathbf{e_3}^T), \tag{8}$$

where $a_1, a_2$, and $a_3$ are coefficients that depend linearly on $k_f$ and the eigenvalues.

Another study by (Painter and Hillen, 2013) has suggested a different derivation of the anisotropic behavior of gliomas from *peanut* or von Mises distribution. In this approach, the macroscopic behavior of glioma cells is derived by considering individual pathways of tumor cells along white matter fibres (Hillen, 2006). For example for the 3D case of von Mises-Fisher distribution the diffusion tensor will have the form of:

$$\mathbf{T} = a_1 I + a_2 \lambda_1 \mathbf{e_1} \mathbf{e_1}^T, \tag{9}$$

where $a_1$ and $a_2$ are two parameters that depend on the sensitivity of cells to pathways, cell turning rate, and the degree of randomized turning. In this paper, we consider Eq. 5 and 7. But the inverse problem formulation of §3 can handle the cases of Eqns. 5, and 8. For the case of Eq. 9, our method requires a slight modification to invert for both $a_1$ and $a_2$ instead of just one of them.

A common model (Hogea et al., 2008b; Konukoglu et al., 2010a,b; Mang et al., 2012; Rekik et al., 2013; Rockne et al., 2009) for the cell mitosis and necrosis in Eq. 1 is the following, self-limiting logistic reaction term:

$$R(c) = \rho c(1 - c), \tag{10}$$

where $\rho$ is the reaction coefficient. This model captures the exponential growth of the tumor cells in the areas of low concentration, and necrosis in areas where $c > 1$. We use $\rho = 2$, in non-dimensional form, to test the method in §5, which corresponds to $\rho = 0.0006$ per day in dimensional form. For the diffusion coefficient we use $k_0 = 0.1$ in non-dimensional form, which corresponds to $k_0 = 0.01 \frac{mm^2}{day}$ in dimensional form [3].

[3] The parameters were selected in the range specified by (Stein et al., 2007)

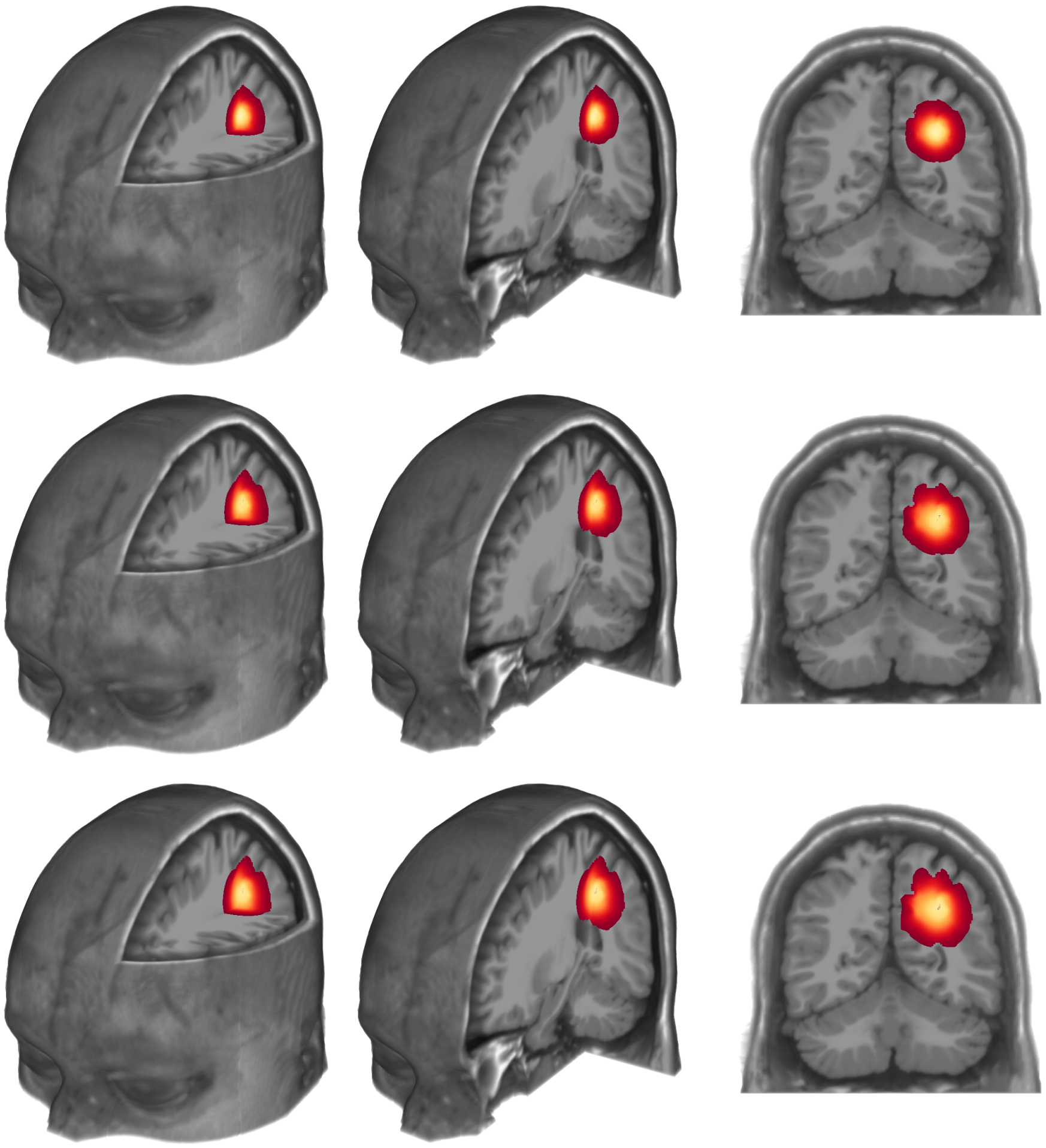

**Fig. 1:** Forward simulation of tumor growth using the reaction-diffusion model of Eq. 1. From top to bottom, the rows show the tumor distribution at t=0,1 and 2, which in dimensional form corresponds to 0, 14, and 28 months, respectively.

As mentioned earlier, with the current clinical imaging technologies only the bulk of the tumor is visible; the full extent of the invasion (physiological tumor boundary) remains undetectable. That is, the tumor concentration is only detectable at locations in which $c^*(\mathbf{x}) > c_d$. Here $c^*(\mathbf{x})$ refers to the spatial tumor distribution and $c_d$ refers to the detection threshold. This threshold depends on the imaging modality used. In (Swanson et al., 2008) the detection threshold was set to $c_d = 0.16$ and $c_d = 0.8$, for T2-weighted (T2w) and T1-weighted contrast enhanced (T1w-Gd) imaging, respectively. However, there is no consensus in the literature as other values such as $c_d = 0.4$ have been used as well (Konukoglu et al., 2010a; Tracqui et al., 1995). A thorough in vivo/vitro experimental study is needed to identify the correct threshold. To the best of our knowledge, no study on this subject has been published yet. Thus, we will consider a range of values for $c_d$ in our tests. In future, if such data becomes available, one can easily substitute its value and run our inversion algorithm. Nothing changes in our formulation, just the value for $c_d$. Based on the detection threshold, we define tumor margin, $m$, to correspond to areas where tumor concentration is below $c_d$, and above a cutoff value of 1%. Because the tumor concentration is continuous we need to use this cutoff to define the margin.

| | |
|---|---|
| $c$ | Normalized tumor concentration |
| $c_i$ | Normalized tumor concentration at $t = i$, $i \in \{0, 1\}$ |
| $c_d$ | Detection threshold |
| $\mathbf{x} = (x, y, z)$ | Spatial location |
| $\mathcal{B}$ | Brain domain |
| $t$ | Time |
| $D$ | Tumor diffusion operator |
| $\mathbf{K}$ | Tumor diffusion tensor |
| $k_0$ | Inhomogeneous diffusion coefficient |
| $k_f$ | Anisotropic diffusion coefficient |
| $\mathbf{T}$ | DTI weighted diffusion tensor |
| $R$ | Tumor reaction operator |
| $\rho$ | Tumor reaction coefficient |
| $d_i$ | Target tumor concentration at $t = i$, $i \in \{0, 1\}$ |
| $O_i$ | Observation operator at $t = i$, $i \in \{0, 1\}$ |
| $p$ | Reconstruction initial condition parametrization |
| $\Phi$ | Reconstruction initial condition parametrization basis function |
| $\beta_p$ | Regularization parameter |
| $\alpha$ | Adjoint variable |

**Table 1:** Basic notations used in this paper.

## 3 Inverse Problem

We consider an inverse problem approach as a modular method of approximating the full extension of tumor invasion, as well as its infiltration rate ($k_f$). The inverse problem is formulated as a PDE constrained optimization problem:

$$\min_{p,k_f} \mathcal{J} := \frac{1}{2}\|O_0 c_0 - d_0\|_2^2 + \frac{1}{2}\|O_1 c_1 - d_1\|_2^2 + \frac{\beta_p}{2}\|p\|_2^2 \tag{11}$$

subject to:

$$\frac{\partial c}{\partial t} - Dc - R(c) = 0 \qquad \text{in U,} \tag{12}$$

$$\frac{\partial c}{\partial n} = 0 \qquad \text{on } \mathit{\Gamma} \times (0,1), \tag{13}$$

$$c_0 - \mathit{\Phi} p = 0 \qquad \text{in } \mathcal{B} \text{ (initial condition).} \tag{14}$$

Here, $O_0$ and $O_1$ are observation operators, $d_1 \in \mathbb{R}^{N_1}$ and $d_0 \in \mathbb{R}^{N_0}$ are the vectors of observation points, $p \in \mathbb{R}^{N_p}$ is a parametrization of tumor distribution at $t = 0$, and $\beta_p$ is a regularization parameter. All the subscripts refer to time (e.g. $O_0$ is observation operator at $t = 0$). The commonly used notations are defined in Table 1.

The observation operators are defined as:

$$O_0 = \begin{cases} 1, & \text{if } \; c_0(\mathbf{x}) \geq c_d \\ 0, & \text{otherwise.} \end{cases}$$

$$O_1 = \begin{cases} 1, & \text{if } \; c_1(\mathbf{x}) \geq c_d \\ 0, & \text{otherwise.} \end{cases}$$

The reason for different observation operators becomes clear now since $c_1$ has a different distribution than $c_0$. As mentioned before, only tumor concentrations above a threshold are observable in medical images. Therefore, the second tumor may have a larger or smaller detectable part.

Moreover, the initial distribution of the tumor in Eq. 14 is parametrized as a superposition of $N_p$ Gaussian distributions. This mapping is done by $\Phi \in \mathbb{R}^{N_0 \times N_p}$.

The corresponding Lagrangian of this problem is given by:

$$L(c, \alpha, p, k_f) = \mathcal{J} + \int_0^1 \int_{\mathcal{B}} \alpha(\frac{\partial c}{\partial t} - Dc - R(c)) d\mathcal{B} dt + \int_{\mathcal{B}} \alpha_0 (c_0 - \Phi p) d\mathcal{B}. \quad (15)$$

The first order optimality conditions can be derived by requiring stationarity of the Lagrangian with respect to state $c$, adjoint $\alpha$ and inversion variables $p,\ k_f$. That is:

$$\frac{\partial L}{\partial c} = 0 \qquad \text{(adjoint equation)} \qquad (16)$$

$$\frac{\partial L}{\partial \alpha} = 0 \qquad \text{(state equation)} \qquad (17)$$

$$\frac{\partial L}{\partial p} = 0 \qquad \text{(inversion equation)} \qquad (18)$$

$$\frac{\partial L}{\partial k_f} = 0 \qquad \text{(inversion equation)} \qquad (19)$$

As a result, the adjoint equations are:

$$-\frac{\partial \alpha}{\partial t} - D\alpha - \frac{\partial R(c)}{\partial c}\bigg|_c \alpha = 0, \qquad \text{in U}, \qquad (20)$$

$$\frac{\partial \alpha}{\partial n} = 0, \qquad \text{on } \Gamma \times (0,1), \qquad (21)$$

$$\alpha_1 + O^T(Oc_1 - d_1) = 0, \qquad \text{in } \mathcal{B} \text{ (initial condition)}. \qquad (22)$$

Similarly, the state equations are:

$$\frac{\partial c}{\partial t} - Dc - R(c) = 0, \qquad \text{in U,} \tag{23}$$

$$\frac{\partial c}{\partial n} = 0, \qquad \text{on } \mathit{\Gamma} \times (0{,}1), \tag{24}$$

$$c_0 - \mathit{\Phi} p = 0, \qquad \text{in } \mathcal{B} \text{ (initial condition).} \tag{25}$$

Finally, the inversion equations are:

$$\beta_p p - \mathit{\Phi}^T \alpha_0 + O_0^T (O_0 \mathit{\Phi} p - d_0) = 0 \tag{26}$$

$$\int_0^1 \int_{\mathcal{B}} (\mathbf{T}\nabla c) \cdot (\nabla \alpha) d\mathcal{B} dt = 0. \tag{27}$$

Equations 20-27 form a system of nonlinear PDEs in space and time. In general, they are not amenable to analytical solutions, and have to be solved with a numerical method. We use a Gauss-Newton method that corresponds to an inexact linearization of the equations. After deriving the optimality conditions, then we discretize and solve the resulting PDEs numerically (Kallivokas et al., 2013). The scheme can be summarized as follows. Assume $c^0$, $\alpha^0$, $p^0$ and $k_f^0$ to be the initial guess of the iterative scheme (throughout the paper superscripts refer to the iteration number of the iterative solver and should not be confused with time, which are always specified by subscripts). The updates to these variables (denoted by $\tilde{c}$, $\tilde{\alpha}$, $\tilde{p}$ and $\tilde{k_f}$) can be found by solving the second order optimality system (i.e. the second variation of the Lagrangian):

$$-\frac{\partial \tilde{\alpha}}{\partial t} - D\tilde{\alpha} - \left.\frac{\partial R}{\partial c}\right|_{c^0} \tilde{\alpha} - \left.\frac{\partial^2 R}{\partial c^2}\right|_{c^0} \tilde{c}\alpha^0 - \nabla \cdot (\tilde{k_f} \mathbf{T} \nabla \alpha^0) = 0 \tag{28}$$

$$\tilde{\alpha}_1 + O^T O \tilde{c}_1 = 0 \tag{29}$$

$$\beta_p \tilde{p} - \mathit{\Phi}^T \tilde{\alpha}_0 + O_0^T (O_0 \mathit{\Phi} \tilde{p} - d_0) = -(\beta_p p^0 - \mathit{\Phi}^T \alpha_0^0 + O_0^T (O_0 \mathit{\Phi} p^0 - d_0)) \tag{30}$$

$$\frac{\partial \tilde{c}}{\partial t} - D\tilde{c} - \left.\frac{\partial R}{\partial c}\right|_{c^0} \tilde{c} - \nabla \cdot (\tilde{k_f} \mathbf{T} \nabla c^0) = 0 \tag{31}$$

$$\tilde{c}_0 - \mathit{\Phi} \tilde{p} = 0 \tag{32}$$

$$\int_0^1 \int_{\mathcal{B}} \big( (\mathbf{T}\nabla \tilde{c}) \cdot (\nabla \alpha^0) + (\mathbf{T}\nabla c^0) \cdot (\nabla \tilde{\alpha}) \big) d\mathcal{B} dt = 0. \tag{33}$$

It is useful to rewrite the previous equations using linear operators (definitions are given in §6):

$$J^T \tilde{\alpha} + N\tilde{c} + Z^T \tilde{k_f} = 0, \tag{34}$$

$$B_p \tilde{p} - \mathit{\Phi}^T \tilde{\alpha} = -g_p, \tag{35}$$

$$J\tilde{c} + W^T \tilde{k_f} - \mathit{\Phi} \tilde{p} = 0, \tag{36}$$

$$Z\tilde{c} + W\tilde{\alpha} = -g_k, \tag{37}$$

which in matrix form can be written as:

$$\begin{bmatrix} N & 0 & J^T & Z^T \\ 0 & B_p & -\Phi^T & 0 \\ J & -\Phi & 0 & W^T \\ Z & 0 & W & 0 \end{bmatrix}_{(c^0,\alpha^0,p^0,k_f^0)} \begin{bmatrix} \tilde{c} \\ \tilde{p} \\ \tilde{\alpha} \\ \tilde{k_f} \end{bmatrix} = \begin{bmatrix} 0 \\ -g_p \\ 0 \\ -g_k \end{bmatrix}_{(c^0,\alpha^0,p^0,k_f^0)} \tag{38}$$

To solve this linear system, we use a reduced Hessian formulation, by eliminating $\tilde{c}$ and $\tilde{\alpha}$ from the system. As a result we obtain:

$$\begin{bmatrix} H_{pp} & H_{pk} \\ H_{kp} & H_{kk} \end{bmatrix} \begin{bmatrix} \tilde{p} \\ \tilde{k_f} \end{bmatrix} = \begin{bmatrix} -g_p \\ -g_k \end{bmatrix}_{(c^0,\alpha^0,p^0,k_f^0)}, \tag{39}$$

where the reduced Hessians (i.e. the individual block matrices on the left hand side in Eq. 39) are defined in §6 (see Eqns. 57-60).

To compute $\tilde{p}$ we solve the following equation using an iterative scheme such as Generalized Minimal Residual (GMRES) or Conjugate Gradient (CG) method. Because the Hessian is symmetric positive definite, we use CG to solve the Schur complement of Eq. 39, given by:

$$(H_{pp} - H_{pk}H_{kk}^{-1}H_{kp})\tilde{p} = H_{pk}H_{kk}^{-1}g_k - g_p. \tag{40}$$

Then $\tilde{k_f}$ can be found by substituting the computed $\tilde{p}$ in the following equation:

$$\tilde{k_f} = -H_{kk}^{-1}H_{kp}\tilde{p} - H_{kk}^{-1}g_k. \tag{41}$$

Since the problem is nonlinear, this update process must be repeated. That is, one should update $(p^0, k_f^0)$ with $(\tilde{p}, \tilde{k_f})$ (using a globalization scheme such as a line search) and repeat the solution process (i.e. solve Eqns. 40 and 41 again).

## 4 Numerical Methods

Following our earlier work (Hogea et al., 2008b), we use a Strang second-order time splitting method (Strang, 1968) to numerically solve the PDEs. Splitting methods provide exclusive benefits for solving complex PDEs that involve different operators, without loss of accuracy. Here we explain how the forward PDE (i.e. Eq. 1) is solved with this method. The other PDEs of the optimality systems in §3 can be solved in a similar fashion. Let $c^n$ denote the tumor distribution at the $n$-th time step. To find $c^{n+1}$ the following steps have to be performed:

1. Solve $\frac{\partial c}{\partial t} = Dc$ over time $\Delta t/2$ using a second order implicit Crank-Nicolson method with $c^n$ as initial condition, to obtain $c^\dagger$.
2. Solve $\frac{\partial c}{\partial t} = R(c)$ over time $\Delta t$ analytically/numerically with $c^\dagger$ as initial condition, to obtain $c^{\dagger\dagger}$.
3. Solve $\frac{\partial c}{\partial t} = Dc$ using a second order implicit Crank-Nicolson method over time $\Delta t/2$ with $c^{\dagger\dagger}$ as initial condition, to obtain $c^{n+1}$.

This scheme can more compactly be written as:

$$c^{n+1} = S_D^{\frac{\Delta t}{2}} S_R^{\Delta t} S_D^{\frac{\Delta t}{2}} c^n, \tag{42}$$

where $S_i^t$, $i = D, R$, is the numerical PDE solver corresponding to each of the above steps. For instance, $S_D^{\frac{\Delta t}{2}} c^n$ corresponds to the first step and applies the diffusion operator to $c^n$ and gives $c^\dagger$ as its output (the interested reader is referred to (Strang, 1968) and (Ropp and Shadid, 2009) for more details). To solve for each of these steps we use a pseudo-spectral method. The space is discretized into $N = 64^3$ nodes and $N_t = 10$ time steps. The space discretization corresponds to a domain of 64 $mm^3$ centered around the tumorous region[4]. Finally, we note that the overall numerical scheme is second order in time and has spectral accuracy in space. We have tested and verified this.

A key advantage of the splitting method is that the nonlinear logistic growth model has an analytical solution. As a result, $S_R^{\Delta t}$ can be computed with an accuracy that is down to machine precision and with negligible computational cost (this not only reduces the computational complexity, but also allows for very accurate approximations to the Hessian operator). For PDEs that involve the derivative of the reaction term, one can simply use a second order Crank-Nicolson method in the absence of an exact solution. For the diffusion part $S_D^{\frac{\Delta t}{2}}$ we use pseudo-spectral method in space and second order Crank-Nicolson scheme in time. For instance, $S_D^{\frac{\Delta t}{2}} c^n$ is solved to obtain $c^\dagger$ as follows:

$$(I - \frac{\Delta t}{4} \nabla \cdot \mathbf{K} \nabla) c^\dagger = (I + \frac{\Delta t}{4} \nabla \cdot \mathbf{K} \nabla) c^n. \tag{43}$$

The implicit scheme of Eq. 43 results in a symmetric system that is solved iteratively using the Conjugate Gradient method. A key advantage of implicit methods is that they are, in contrast to explicit methods, unconditionally stable. On the downside, one must solve a system of equations at each iteration.

*Hessian Preconditioner.* Here, we briefly explain how we solve and precondition Eq. 40. Solving this equation is one of the major bottlenecks of the algorithm, as each application of the Hessian involves the solution of the forward and adjoint PDEs in time. To prevent any ambiguity, we refer to $(H_{pp} - H_{pk} H_{kk}^{-1} H_{kp})\tilde{p}$ as the *Hessian matvec* for a given $\tilde{p}$.

We precondition Eq. 40 to increase its convergence rate. As a result, the number of necessary Hessian matvecs will decrease, thus reducing time to solution. A good preconditioner is one that has low overhead of forming and applying, and yet increases the convergence rate considerably. An extreme case for a preconditioner is the inverse of the Hessian operator, which would decrease the number of iterations to one, but is too expensive to form and apply. On the other hand, the identity operator is a preconditioner with no overhead and obviously no effect. With this in mind, we preconditioned Eq. 40 with an analytical approximation to $H_{pp}$. The analytic expression is derived by approximating numerical operators of $H_{pp}$ (defined in Eq. 57). For this, we use the average value of the variable diffusion coefficient tensor $\mathbf{K}(\mathbf{x})$. This allows an analytical solution to the diffusion operator in the splitting scheme. The other terms such as the reaction operator are solved numerically as before. The fact that we are not solving an implicit system of equations for the diffusion operator reduces the overhead of forming and applying this preconditioner considerably. This is the only approximation that we are doing and the rest of the operators are the same as the numerical Hessian. This preconditioner performs very well and reduces the number of CG iterations to solve Eq. 57 by a factor of 5 (from about 25 to 5).

[4] This area can be increased for highly infiltrative tumors that spread through a larger portion of the brain.

Since we are solving a nonlinear optimization problem, Eq. 40 has to be repeatedly solved (we refer to these repetitions as Newton iterations). The number of Newton iterations depends on how far the initial guess of the inversion parameters are from the actual (true) solution. The fact that the Hessian preconditioner performs very well, shows that it is actually a good approximation to the true Hessian. As a result, one can use it to first solve for a good starting point with negligible computational overhead. This means that one first solves Eq. 40 and Eq. 41 to compute $p^0$ and $k_f^0$ by using analytical approximations to the Hessian operators, instead of using their true numerical expressions. The approximation is exactly the same as explained above; that is, we only approximate the costly diffusion operator with an analytical expression by using the average value of the diffusion tensor. With these starting values of $p^0$ and $k_f^0$ we then solve Eq. 40 and Eq. 41 with the correct numerical forms of the Hessian operators. This process is quite effective as it reduces the number of Newton iterations by a factor of 3 (from about 15 to 5).

*Regularization Parameter.* Inherent in any inverse problem, is a regularization operator that limits noise amplification. We use a Tikhonov-type regularization model in Eq. 11 for this purpose (Tikhonov, 1963). If the problem is under-regularized (i.e. $\beta_p \ll 1$), the solution will be dominated by the inherent noise in the data. On the other hand, over-regularization (i.e. $\beta_p \gg 1$) will result in a poor fit to the data. Therefore, it is necessary to choose an optimal regularization parameter, somewhere in-between these two extremes. It is possible to choose $\beta_p$ heuristically. However, we will use a systematic approach to compute an optimal $\beta_p$. The L-curve method is one of the widely used methods for this purpose (Hansen, 1992, 1999). We provide a plot of an exemplary L-curve in Fig. 2. We obtained $\beta_p = 0.01$ as an optimal regularization parameter for test case 2, with $c_d = 0.2$, and $\eta = 5\%$ (the test case will be defined in §5). This process needs to be repeated for new datasets. Although this method has an expensive computational cost, it allows a systematic selection of the regularization parameter. It is also possible to do this process inside the iterative solver, where one starts from the tail of the L-Curve (large regularization) and reduces the regularization parameter as the solution con-verges (Fathi et al., 2015). In this fashion, part of the computational cost will be overlapped in the non-linear iterations.

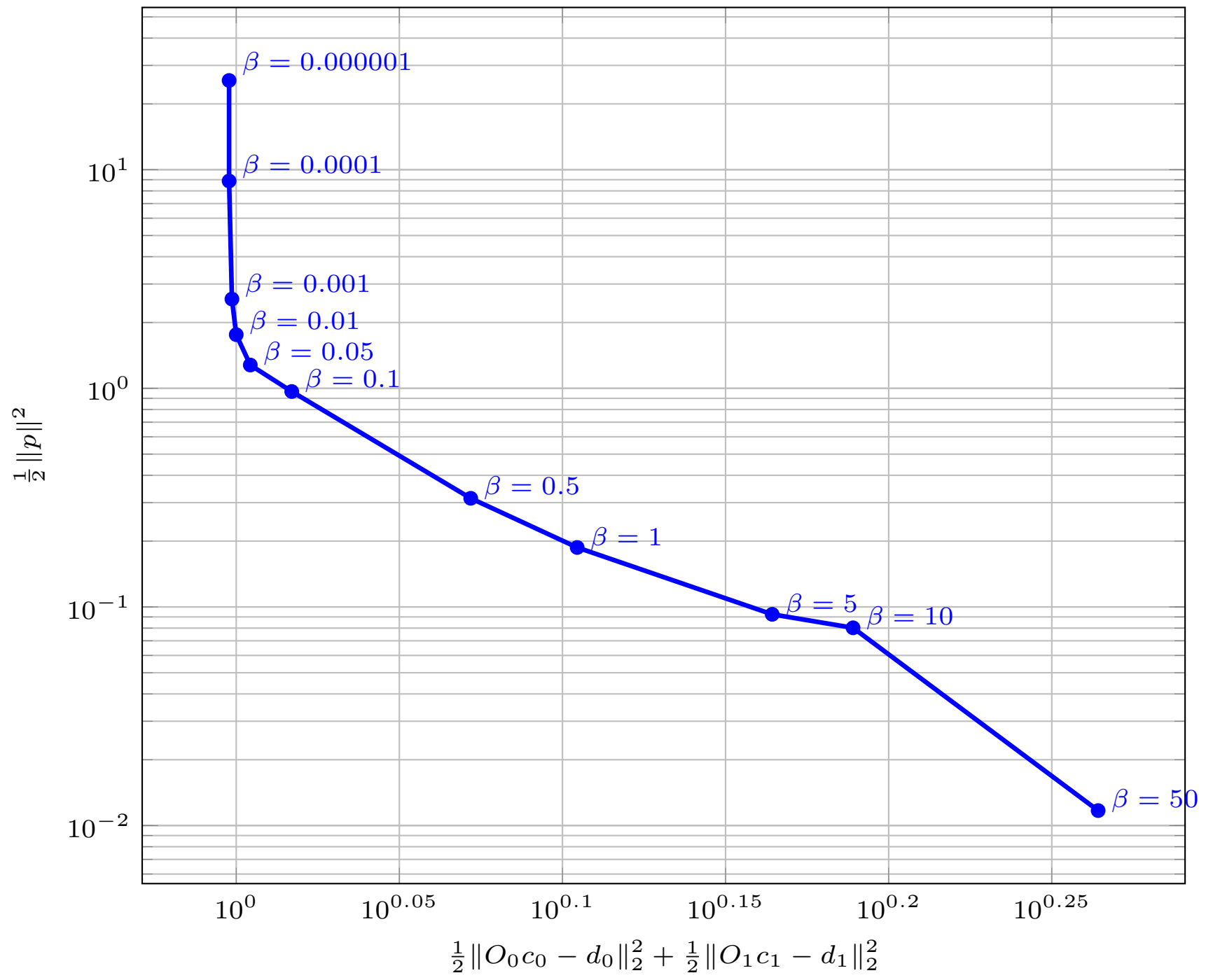

**Fig. 2:** L-curve for choosing the regularization parameter. The corner of the L-curve yields $\beta_p = 0.01$ as the optimal regularization parameter for $c_d = 0.2$ and $\eta = 5\%$.

## 5 Results

*Setup.* We consider synthetically created tumor distributions as target. We apply our inversion algorithm to noisy, partial observations of these targets, limited to two time frames, and compute the corresponding reconstructions. This synthetic analysis allows us to compute the exact errors between the reconstructions and the target, which would otherwise not be possible. As mentioned before, only part of the target tumor is detectable. Our observation operator (i.e. $O_0$ and $O_1$) captures this, by selecting tumors at points that have a concentration above the detection threshold. Since there is no consensus in the literature about the exact value of this detection threshold (Konukoglu et al., 2010a; Swanson et al., 2008; Tracqui et al., 1995), we will consider different values to test our algorithm. Moreover, the data derived from any imaging modality will contain noise. To account for this, we add different levels of white noise (denoted by $\eta$) to the observed data (target distribution).

**Table 2:** Reconstruction results for test case 1. For this test case, the target value of the anisotropic diffusion coefficient $k_f$ is known. The reconstruction errors, $\epsilon$ , are given at different time frames, for different observation thresholds $c_d$ and different noise levels $\eta$.

| $c_d$ | $\eta$ | $\epsilon_0$ | $JI_0$ | $\epsilon_1$ | $JI_1$ | $\epsilon_2$ | $JI_2$ |
|---|---|---|---|---|---|---|---|
| 0.10 | 1% | 5.7e-02 | 0.794 | 3.5e-02 | 0.840 | 3.4e-02 | 0.821 |
| 0.10 | 5% | 6.2e-02 | 0.806 | 3.9e-02 | 0.854 | 3.6e-02 | 0.835 |
| 0.10 | 10% | 7.5e-02 | 0.818 | 5.4e-02 | 0.865 | 4.9e-02 | 0.846 |
| 0.20 | 1% | 6.1e-02 | 0.816 | 4.2e-02 | 0.857 | 4.1e-02 | 0.840 |
| 0.20 | 5% | 6.7e-02 | 0.825 | 4.7e-02 | 0.866 | 4.3e-02 | 0.849 |
| 0.20 | 10% | 8.3e-02 | 0.832 | 6.3e-02 | 0.874 | 5.7e-02 | 0.857 |
| 0.40 | 1% | 6.8e-02 | 0.845 | 4.4e-02 | 0.889 | 4.0e-02 | 0.871 |
| 0.40 | 5% | 7.7e-02 | 0.849 | 5.3e-02 | 0.896 | 4.6e-02 | 0.877 |
| 0.40 | 10% | 1.0e-01 | 0.852 | 7.7e-02 | 0.899 | 6.8e-02 | 0.880 |

To match our implementation to virtual brain anatomy, we use the BrainWeb atlas (spatial resolution: $1mm \times 1mm \times 1mm$) (Cocosco et al., 1997). Moreover, we use the diffusion tensor imaging data provided by the LONI lab of the University of Southern California (Mori et al., 2008).

To test our inversion algorithm, we consider the following four test cases:

Test case 1: Reconstruct full tumor distribution with known anisotropic diffusion coefficient $k_f$ (using Eq. 5 for the diffusion tensor).
Test case 2: Reconstruct full tumor distribution as well as anisotropic diffusion coefficient $k_f$ (using Eq. 5 for the diffusion tensor).
Test case 3: Same as test case 2, but for a multifocal tumor (using Eq. 5 for the diffusion tensor).
Test case 4: Reconstruct full tumor distribution as well as anisotropic diffusion coefficient $k_f$ (using Eq. 7 for the diffusion tensor).

The inversion in all the test cases is derived by partially observed, noisy data of the target tumor distribution at two consecutive time frames of $t = 0$ and $t = 1$, respectively. Test case 1 is rather unrealistic, since there is currently no way to measure $k_f$ in vivo. However, in test case 2 we assume no knowledge of the anisotropic diffusion rate. We include it as an unknown inversion parameter in our algorithm. Other researchers have considered manual tuning of $k_f$ and then selected the one that has the best fit (Jbabdi et al., 2005; Sodt et al., 2014). However, our approach is a systematic one and does not need any manual parameter tuning. In the last test case we consider a multifocal tumor. This test case is to demonstrate that our inversion algorithm works irrespective of mono- or multifocality of the tumor. Moreover, to show that our method works with other diffusion tensors, we test the inversion using Eq. 7 for **T**.

*Performance Measures.* To quantify the reconstruction performance, we report the mismatch between the target and reconstructed tumor distributions:

$$\epsilon_i = \frac{\|c_i - c_i^*\|_2}{\|c_i^*\|_2}, \tag{44}$$

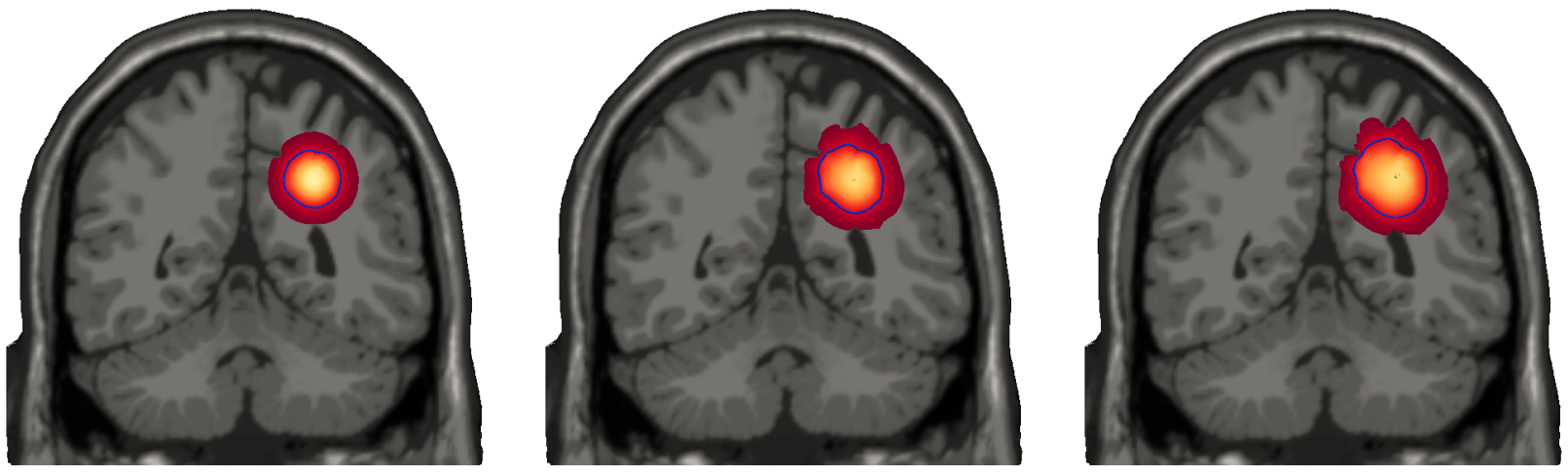

**(a)** Target Tumor at t=0 **(b)** Target Tumor at t=1 **(c)** Target Tumor at t=2

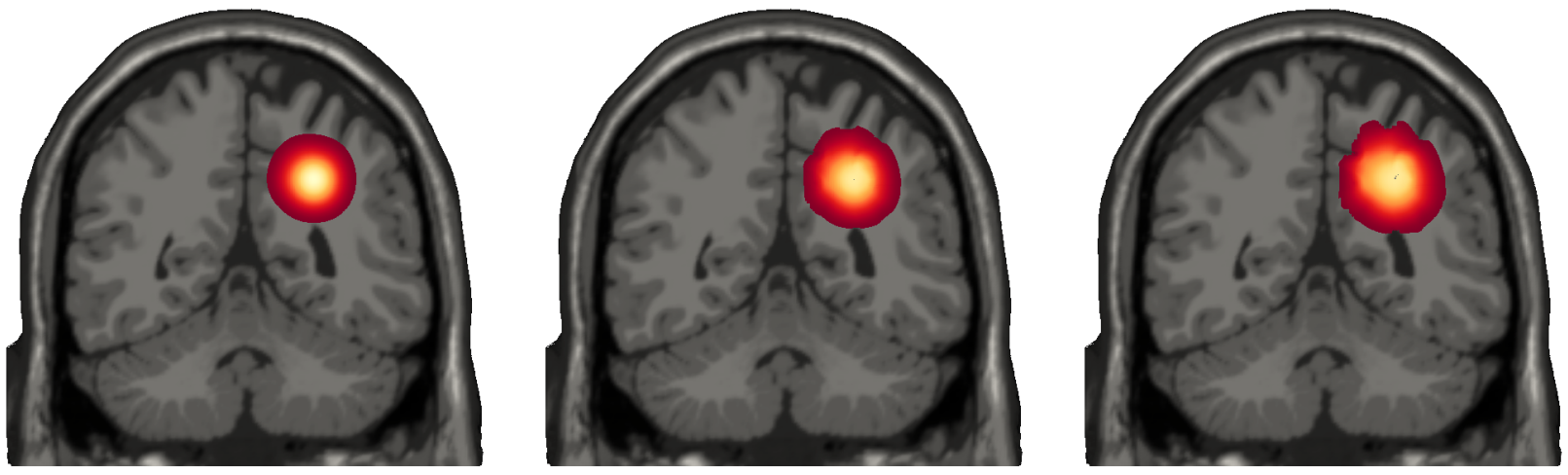

**(d)** Reconstructed Tumor at t=0 **(e)** Reconstructed Tumor at t=1 **(f)** Reconstructed Tumor at t=2

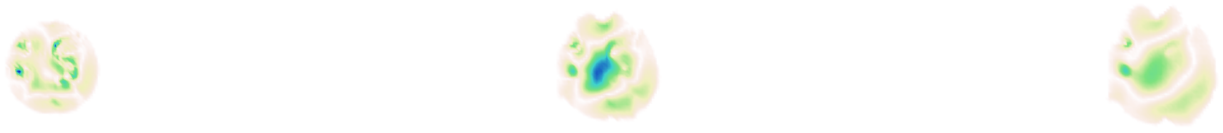

**(g)** Rescaled Difference at t=0 **(h)** Rescaled Difference at t=1 **(i)** Rescaled Difference at t=2

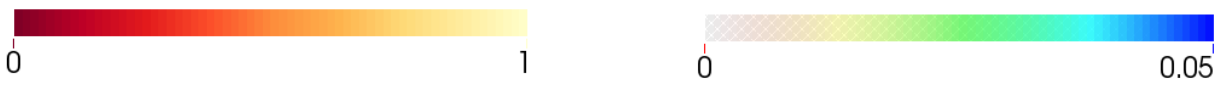

**(j)** **(k)**

**Fig. 3:** Reconstruction results for test case 2. The top row shows the target tumor distribution at different points in time. The three columns show the tumor distribution at $t = 0, 1$, and 2 which in dimensional form corresponds to 0, 14, and 28 months, respectively. The blue contour indicates an observable tumor concentration of $c_d = 0.2$. The reconstruction is driven by noisy observations ($\eta = 5\%$) of the cell density on the inside of this contour (at $t = 0$ and $t = 1$). The second row shows the reconstruction results at the corresponding times. The reconstruction relative errors are 6.6%, 4.5% and 5.3%, respectively. The last row shows the same difference but with a rescaled colormap (k) of the original colormap (j).

**Table 3:** Detailed reconstruction results for test case 2. The reconstruction errors, $\epsilon$ , are given at different time frames, for different observation thresholds $c_d$ and different noise levels $\eta$.

| $c_d$ | $\eta$ | $\epsilon_{kf}$ | $\epsilon_0$ | $JI_0$ | $\epsilon_1$ | $JI_1$ | $\epsilon_2$ | $JI_2$ |
|---|---|---|---|---|---|---|---|---|
| 0.10 | 1% | 5.4e-02 | 5.7e-02 | 0.793 | 3.4e-02 | 0.804 | 4.6e-02 | 0.800 |
| 0.10 | 5% | 5.1e-02 | 6.2e-02 | 0.806 | 3.8e-02 | 0.818 | 4.8e-02 | 0.814 |
| 0.10 | 10% | 6.8e-02 | 7.5e-02 | 0.818 | 5.3e-02 | 0.831 | 5.8e-02 | 0.827 |
| 0.20 | 1% | 1.3e-02 | 6.1e-02 | 0.816 | 4.2e-02 | 0.829 | 5.2e-02 | 0.824 |
| 0.20 | 5% | 1.3e-02 | 6.7e-02 | 0.825 | 4.7e-02 | 0.840 | 5.4e-02 | 0.834 |
| 0.20 | 10% | 1.3e-02 | 8.3e-02 | 0.832 | 6.3e-02 | 0.848 | 6.5e-02 | 0.843 |
| 0.30 | 1% | 2.4e-03 | 6.4e-02 | 0.828 | 4.5e-02 | 0.844 | 5.3e-02 | 0.837 |
| 0.30 | 5% | 3.3e-03 | 7.1e-02 | 0.834 | 5.0e-02 | 0.853 | 5.6e-02 | 0.845 |
| 0.30 | 10% | 1.2e-02 | 9.1e-02 | 0.839 | 6.9e-02 | 0.862 | 6.9e-02 | 0.853 |
| 0.40 | 1% | 1.4e-02 | 6.8e-02 | 0.844 | 4.4e-02 | 0.866 | 5.2e-02 | 0.858 |
| 0.40 | 5% | 3.6e-02 | 7.6e-02 | 0.849 | 5.2e-02 | 0.873 | 5.5e-02 | 0.864 |
| 0.40 | 10% | 6.6e-02 | 9.8e-02 | 0.853 | 7.5e-02 | 0.879 | 7.2e-02 | 0.869 |

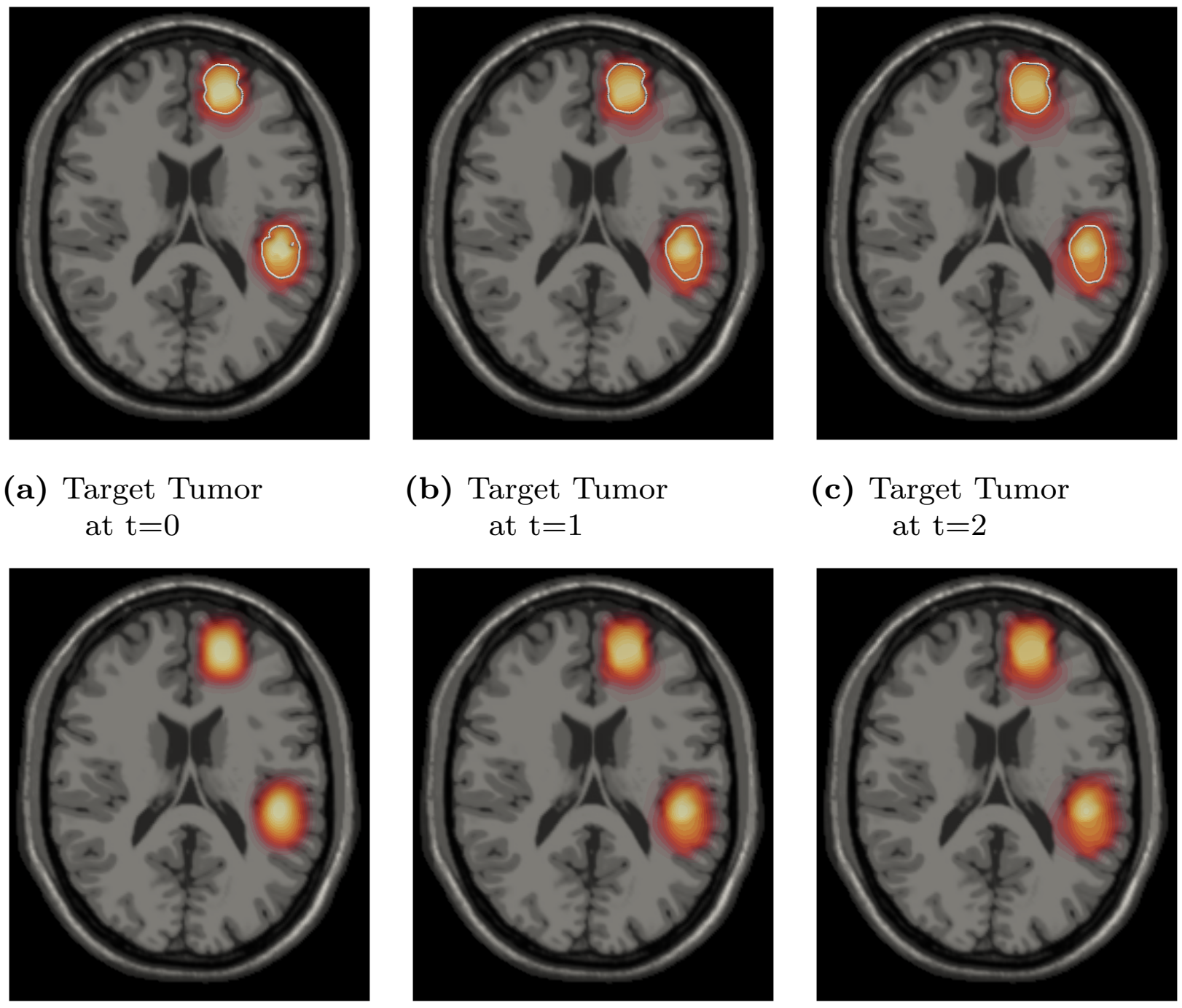

**(a)** Target Tumor at t=0 **(b)** Target Tumor at t=1 **(c)** Target Tumor at t=2

**(d)** Reconstructed Tumor at t=0 **(e)** Reconstructed Tumor at t=1 **(f)** Reconstructed Tumor at t=2

**Fig. 4:** Reconstruction results for multifocal test case 3. The top row shows the target tumor distribution at different times. The three columns show the tumor distribution at $t = 0, 1$, and 2 which in dimensional form corresponds to 0, 5, and 10 months, respectively. The white contour indicates the observable tumor concentration of $c_d = 0.2$. The bottom row shows the reconstruction results at $t = 0, 1$, and 2 for $\eta = 5\%$. The reconstruction errors are 10.3%, 6.97% and 6.56%, respectively.

**Table 4:** Detailed reconstruction results for multifocal test case 3. The reconstruction errors, $\epsilon$, are given at different time frames, for different observation thresholds $c_d$ and different noise levels $\eta$.

| $c_d$ | $\eta$ | $\epsilon_{kf}$ | $\epsilon_0$ | $JI_0$ | $\epsilon_1$ | $JI_1$ | $\epsilon_2$ | $JI_2$ |
|---|---|---|---|---|---|---|---|---|
| 0.10 | 1% | 1.0e-02 | 9.7e-02 | 0.695 | 6.1e-02 | 0.748 | 5.5e-02 | 0.726 |
| 0.10 | 5% | 2.3e-02 | 1.0e-01 | 0.695 | 6.6e-02 | 0.753 | 6.1e-02 | 0.726 |
| 0.10 | 10% | 3.8e-02 | 1.1e-01 | 0.688 | 7.9e-02 | 0.747 | 7.4e-02 | 0.722 |
| 0.20 | 1% | 8.0e-02 | 1.0e-01 | 0.748 | 7.2e-02 | 0.795 | 6.6e-02 | 0.775 |
| 0.20 | 5% | 8.0e-02 | 1.1e-01 | 0.742 | 7.8e-02 | 0.790 | 7.3e-02 | 0.771 |
| 0.20 | 10% | 8.1e-02 | 1.2e-01 | 0.739 | 9.2e-02 | 0.790 | 8.7e-02 | 0.769 |
| 0.30 | 1% | 7.4e-02 | 1.2e-01 | 0.764 | 9.0e-02 | 0.801 | 8.5e-02 | 0.785 |
| 0.30 | 5% | 6.5e-02 | 1.2e-01 | 0.766 | 9.6e-02 | 0.804 | 9.1e-02 | 0.788 |
| 0.30 | 10% | 5.5e-02 | 1.4e-01 | 0.760 | 1.1e-01 | 0.804 | 1.1e-01 | 0.786 |
| 0.40 | 1% | 1.7e-01 | 1.3e-01 | 0.763 | 1.0e-01 | 0.810 | 9.4e-02 | 0.790 |
| 0.40 | 5% | 1.7e-01 | 1.3e-01 | 0.766 | 1.0e-01 | 0.807 | 9.7e-02 | 0.789 |
| 0.40 | 10% | 1.7e-01 | 1.4e-01 | 0.773 | 1.1e-01 | 0.811 | 1.1e-01 | 0.793 |

**Table 5:** Detailed reconstruction results for multifocal test case 4. This test case is different, in that we use Eq. 7 for **T**. The reconstruction errors, $\epsilon$, are given at different time frames, for different observation thresholds $c_d$ and different noise levels $\eta$.

| $c_d$ | $\eta$ | $\epsilon_{kf}$ | $\epsilon_0$ | $JI_0$ | $\epsilon_1$ | $JI_1$ | $\epsilon_2$ | $JI_2$ |
|---|---|---|---|---|---|---|---|---|
| 0.20 | 1% | 1.09e-02 | 7.58e-02 | 0.806 | 4.81e-02 | 0.820 | 5.69e-02 | 0.814 |
| 0.20 | 5% | 2.39e-02 | 8.36e-02 | 0.815 | 5.51e-02 | 0.829 | 6.07e-02 | 0.824 |
| 0.20 | 10% | 4.23e-02 | 1.04e-01 | 0.823 | 7.57e-02 | 0.839 | 7.49e-02 | 0.833 |
| 0.40 | 1% | 6.66e-03 | 6.96e-02 | 0.834 | 4.60e-02 | 0.857 | 5.74e-02 | 0.848 |
| 0.40 | 5% | 8.65e-03 | 7.52e-02 | 0.838 | 5.03e-02 | 0.863 | 5.90e-02 | 0.853 |
| 0.40 | 10% | 1.01e-02 | 9.02e-02 | 0.841 | 6.60e-02 | 0.867 | 6.98e-02 | 0.857 |

where $c_i^*$ is the target tumor distribution and $c_i$ is the reconstructed one, at time $t = i$. This mismatch error is reported for the two consecutive time frames, at which the data was observed (i.e. $t = 0$ and $t = 1$). We also report the error for an additional time point $t = 2$. This test is included to study how well the reconstruction captures the growth pattern of an untreated tumor (i.e. the tumor growth prediction capabilities of our method). Furthermore, we provide values for the Jaccard Index (JI) to assess how well our method approximates the tumor margin:

$$\mathrm{JI}_i = \frac{|H(m_i) \cap H(m_i^*)|}{|H(m_i) \cup H(m_i^*)|}, \tag{45}$$

where $H$ is the Heaviside function, $m_i$ is the reconstructed margin at time $t = i$, and $m_i^*$ is the target one. Note that we are using the tumor margin, $m$, instead of the full tumor distribution, $c$. This is because the tumor margin does not involve areas above the detection threshold. Those areas have a JI of one, and including them will increase JI artificially. This is specially important for cases where $c_d$ is small, and using $m$ will give a better metric of how well the margin is reconstructed.

For other test cases, we also report the error in the reconstruction of $k_f$ compared to the target anisotropy coefficient $k_f^*$:

$$\epsilon_{kf} = \frac{|k_f - k_f^*|}{|k_f^*|}, \tag{46}$$

The target tumor distribution $c^*$ is shown in Fig. 1 at these three time frames. The results are computed by solving the forward problem in 3D. For the test case 1, 2, and 4 we perform the reconstruction in 3D. For the test case 3 (multifocal tumor) we limit the reconstruction to a 2D slice.

*Numerical Results.* The reconstruction results for test case 1 are shown in Table 2. We report the relative reconstruction error (computed according to Eq. 44) at different time frames, for different threshold values, $c_d$, and different noise levels, $\eta$. Note that only the two time frames $t = 0, 1$ are used to solve the inverse problem. The error given at $t = 2$ demonstrates how well the reconstruction matches the target if the tumor was left untreated. The reconstruction error increases with an increasing detection threshold. This behavior is expected, since a higher threshold means less data is available to drive the inversion. Moreover, the reconstruction error is higher, as noise is increased. However, the difference is not significant, because the regularization operator counter balances noise amplifications.

But what if the value of the anisotropic diffusion rate is unknown as well? This is exactly what we study in test case 2. The corresponding errors are shown in Table 3. As one can see, the inversion algorithm approximates the correct value of $k_f$ with a very good accuracy. As a result, the errors in tumor distribution reconstruction at $t = 0$ and $t = 1$ are comparable to those reported for test case 1. Figure 3 illustrates how well the reconstruction compares to the target distribution for test case 2 with $c_d = 0.20$ and $\eta = 5\%$. The target distributions are shown with a blue contour marking the area where the tumor concentration is equal to the detection threshold. The reconstruction, shown in the second row, is computed using only the target cell density that is within this contour. As can be seen, the reconstruction qualitatively captures the growth pattern very well.

As a proof-of-concept, we test how the inversion performs when we consider a multifocal tumor as shown in Fig. 4. This test case is performed in 2D, in contrast to case 1 and 2. The corresponding errors are reported in Table 4. Overall, we observe the same behavior as in the previous test cases. To show that our scheme works for other diffusion tensors, we consider the case of Eq. 7 for $\mathbf{T}$. In this case, only the dominant eigendirection is considered. The results are given in Table 5.

## 6 Conclusions

We presented an inverse problem formulation to determine the full extent of tumor infiltration in the brain based on a PDE-constrained optimization problem. The key quantities of interest are $(i)$ the full extent of tumor invasion, and $(ii)$ the rate of anisotropic diffusion. We used a nonlinear reaction-diffusion model for glioma growth, and solved the optimization problem with a reduced space Hessian method. State of the art numerical techniques were presented to speed up the time to solution. The design criteria for these techniques were low computational cost and robustness. We tested the resulting algorithm, using synthetic tumors with different levels of noise and different detection thresholds.

We want to emphasize that these are preliminary results and a significant amount of work is still necessary, before this scheme can be applied to real datasets. For instance, the use of more complex tumor growth models has to be investigated (Habib et al., 2003; Hawkins-Daarud et al., 2013; Hogea et al., 2008b; Lima et al., 2014), accounting e.g. for mass effect, edema, necrosis, angiogenesis, chemotaxis, and haptotaxis. More importantly, one has to find ways to experimentally verify the forward model as well as the reconstruction results. Moreover, effects of treatment (chemo- and/or radiotherapy) should be added to the model, as serial scans of untreated human subjects are rare.

## Appendix A: Operator Definitions

The definitions of operators in Eqs 34, 35, 36, and 37 are as follows:

$$J^T\tilde{\alpha} := -\frac{\partial\tilde{\alpha}}{\partial t} - D\tilde{\alpha} - \left.\frac{\partial R}{\partial c}\right|_{c^0}\tilde{\alpha} + \int_0^1 \delta(t-T)\tilde{\alpha}dt \tag{47}$$

$$N\tilde{c} := -\left.\frac{\partial^2 R}{\partial c^2}\right|_{c^0}\tilde{c}\alpha^0 + O^T O\int_0^1 \delta(t-T)\tilde{c}dt \tag{48}$$

$$Z^T\tilde{k_f} := -\nabla\cdot(\tilde{k_f}\mathbf{T}\nabla\alpha^0) \tag{49}$$

$$g_p := Bu - \Phi^T\alpha_0^0 \tag{50}$$

$$Bu := (\beta_p + O_0^T O_0\Phi)p \tag{51}$$

$$Z\tilde{c} := \int_0^1\int_\Omega (\mathbf{T}\nabla\tilde{c})\cdot(\nabla\alpha^0)d\Omega dt \tag{52}$$

$$W\tilde{\alpha} := \int_0^1\int_\Omega (\mathbf{T}\nabla c^0)\cdot(\nabla\tilde{\alpha})d\Omega dt \tag{53}$$

$$g_k := \int_0^1\int_\Omega (\mathbf{T}\nabla c^0)\cdot(\nabla\alpha^0)d\Omega dt \tag{54}$$

$$J\tilde{c} := \frac{\partial\tilde{c}}{\partial t} - D\tilde{c} - \left.\frac{\partial R}{\partial c}\right|_{c^0}\tilde{c} + \int_0^1 \delta(t)\tilde{c}dt \tag{55}$$

$$W^T\tilde{k_f} := -\nabla\cdot(\tilde{k_f}\mathbf{T}\nabla c^0) \tag{56}$$

The reduced Hessians of Eq. 39 are defined as follows:

$$H_{pp} = B + \Phi^T J^{-T} N J^{-1} \Phi, \tag{57}$$
$$H_{pk} = -\Phi^T J^{-T}(-Z^T + N J^{-1} W^T), \tag{58}$$
$$H_{kp} = (Z - W J^{-T} N) J^{-1} \Phi, \tag{59}$$
$$H_{kk} = -Z J^{-1} W^T + W J^{-T}(-Z^T + N J^{-1} W^T). \tag{60}$$

To compute the matvec of $H_{pp}\tilde{p}$ one needs to take the following steps:

1. $J^{-1}\Phi$:
   Solve Eq. 31 with $\tilde{k_f} = 0$ and initial condition of $\tilde{c}_0 = \Phi\tilde{p}$ to get $\tilde{c}$

2. $-J^{-T} N J^{-1} \Phi$:
   Solve Eq. 28 with $\tilde{k_f} = 0$ and initial condition of $\tilde{\alpha}_1 = -O^T O \tilde{c}_1$ to get $\tilde{\alpha}$

3. $B + \Phi^T J^{-T} N J^{-1} \Phi$:
   Compute $-\Phi^T \alpha_0$ and add $Bu$

To compute the matvec of $H_{pk}\tilde{k_f}$ one needs to take the following steps:

1. $J^{-1} W^T$:
   Solve Eq. 31 with zero initial condition to get $\tilde{c}$

2. $J^{-T}(-Z^T + N J^{-1} W^T)$:
   Solve Eq. 28 with initial condition of $\tilde{\alpha}_1 = -O^T O \tilde{c}_1$ to get $\tilde{\alpha}$

3. $-\Phi^T J^{-T}(-Z^T + N J^{-1} W^T)$:
   Compute $-\Phi^T \alpha_0$

To compute the matvec of $H_{ku}\tilde{p}$ one needs to take the following steps:

1. $J^{-1}\Phi$:
   Solve Eq. 31 with $\tilde{k_f} = 0$ and initial condition of $\tilde{c}_0 = \Phi\tilde{p}$ to get $\tilde{c}$

2. $Z J^{-1} \Phi$:
   Compute $\int_0^1 \int_\Omega (\mathbf{T}\nabla\tilde{c}) \cdot (\nabla\alpha^0) d\Omega dt$

3. $-J^{-T} N J^{-1} \Phi$:
   Solve Eq. 28 with $\tilde{k_f} = 0$ and initial condition of $\tilde{\alpha}_1 = -O^T O \tilde{c}_1$ to get $\tilde{\alpha}$

4. $-W J^{-T} N J^{-1} \Phi$:
   Compute $\int_0^1 \int_\Omega (\mathbf{T}\nabla c^0) \cdot (\nabla\tilde{\alpha}) d\Omega dt$

5. Add 2 and 4

To compute the matvec of $H_{kk}\tilde{k_f}$ one needs to take the following steps:

1. $J^{-1}W^T$:
   Solve Eq. 31 with zero initial condition to get $\tilde{c}$

2. $ZJ^{-1}W^T$:
   Compute $\int_0^1 \int_\Omega (\mathbf{T}\nabla\tilde{c}) \cdot (\nabla\alpha^0) d\Omega dt$

3. $J^{-T}(-Z^T + NJ^{-1}W^T)$:
   Solve Eq. 28 with initial condition of $\tilde{\alpha}_1 = -O^T O\tilde{c}_1$ to get $\tilde{\alpha}$

4. $WJ^{-T}(-Z^T + NJ^{-1}W^T)$:
   Compute $\int_0^1 \int_\Omega (\mathbf{T}\nabla c^0) \cdot (\nabla\tilde{\alpha}) d\Omega dt$

5. Add 2 and 4.

## Appendix B: Fictitious Domain Method

We use a fictitious domain method in which the original brain domain, $\mathcal{B}$ , is extended to a cubic box, denoted by $\Omega$ (Hogea et al., 2008b; Mang et al., 2012; Tracqui et al., 1995). The original homogeneous boundary conditions imposed on $\Gamma$ can be satisfied using a penalty method (Del Pino and Pironneau, 2003). To do so we define a new diffusion coefficient $\mathbf{K}_\epsilon(\mathbf{x})$, $\mathbf{x} \in \Omega$ as follows:

$$\mathbf{K}_\epsilon(\mathbf{x}) = \begin{cases} \mathbf{K}(\mathbf{x}), & \text{if } \ \mathbf{x} \in \mathcal{B} \\ \epsilon\ \mathbf{K}(\mathbf{x}), & \text{otherwise (i.e. } \mathbf{x} \in \Omega \backslash \bar{\mathcal{B}}) \end{cases}$$

where the penalty parameter $\epsilon$, is a small positive number. The actual boundary condition on $\Gamma$ will be satisfied in the limit of $\epsilon \to 0$ (Del Pino and Pironneau, 2003). The original boundary conditions can be re-imposed on the extended cubic domain, $\Omega$, for both the forward Eq. 1 and adjoint equation Eq. 20.

## Acknowledgments

We would like to thank Thomas Hillen for the helpful discussion on the anisotropic diffusion of gliomas. We would like to also thank Florian Tramnitzke for contributing to this work during his internship in our group.